\documentclass[11pt]{article}

\usepackage{amsmath}
\usepackage{amssymb}
\usepackage{amsthm}
\usepackage{indentfirst}
\usepackage{cite}

\newtheorem{thm}{Theorem}[section]
\newtheorem{lem}{Lemma}[section]
\newtheorem{cor}{Corollary}[section]
\newtheorem{rem}{Remark}[section]
\newtheorem{exm}{Example}

\numberwithin{equation}{section}

\title{\bf ASSOCIATED DERIVED INVARIANTS\\ of Geometric Mappings}
\author{Nenad O. Vesi\'c\footnote{Faculty of Science and Mathematics, Ni\v s, Serbia,
Serbian Ministry of Education, Science and Technological
Development, Grant No. 174012}}
\date{}

\makeatletter
\def\maketag@@@#1{\hbox{\m@th\normalfont\normalsize#1}}
\makeatother

\begin{document}
  \maketitle

  \begin{abstract}
    General invariants of a geometric mapping of a symmetric affine
    connection space are obtained in this paper. These invariants
    are generalizations of the previous obtained basic invariants
    (see \cite{vesicgeninv1}). Moreover, these invariants are
    related with the Thomas projective parameter and the Weyl
    projective tensor.\\[2pt]


    \noindent\textbf{$2010$ Math. Subj. Classification:} 53A55,
    53B05, 53C15
  \end{abstract}

  \section{Introduction}

  This paper is devoted to generalize the Thomas
  projective parameter \cite{thomas1} and the Weyl projective tensor
  \cite{weyl} as invariants of a mapping of an affine connection
  space. These invariants are primary generalized in
  \cite{vesicgeninv1} but it is obtained these invariants are not
  unique ones in this paper. Therefore, we are interested to obtain
  some invariants of geometric mappings different of the invariants
  from \cite{vesicgeninv1} in this paper.

  \subsection{Spaces of affine connection}

  An $N$-dimensional manifold $\mathcal M^N$ equipped with an affine
  connection $\nabla$ with torsion is called the non-symmetric affine
  connection space $\mathbb{GA}_N$. The affine connection
  coefficients of the affine connection $\nabla$ are $L^i_{jk}$ and
  they are non-symmetric by indices $j$ and $k$.

  The geometrical object

  \begin{equation}
        L^i_{\underline{jk}}=\frac12\big(L^i_{jk}+L^i_{kj}\big),
  \end{equation}

  \noindent satisfies the equation

  \begin{equation*}
    L^{i'}_{\underline{j'k'}}=x^{i'}_ix^j_{j'}x^k_{k'}L^i_{\underline{jk}}+
    x^{i'}_{i}x^i_{j'k'},
  \end{equation*}

  \noindent so it is the affine connection coefficient of a
  symmetric affine connection space $\mathbb A_N$.

    The manifold $\mathcal M^N$ equipped with an affine connection
  $\widetilde\nabla$ whose coefficients are $L^i_{\underline{jk}}$ is
  \textbf{the associated space $\mathbb A_N$} (of the space
  $\mathbb{GA}_N$).

  There is one kind of covariant derivation with regard to a
  symmetric affine connection:

    \begin{equation}
    a^i_{j|k}=a^i_{j,k}+L^i_{\underline{\alpha k}}a^\alpha_j-
    L^\alpha_{\underline{jk}}a^i_\alpha,
    \label{eq:covderivativesim}
  \end{equation}

  \noindent for a tensor $a^i_j$ of the type $(1,1)$ and the partial
  derivative $\partial/\partial x^i$ denoted by comma.

  From the corresponding Ricci-type identity, it is obtained one
  curvature tensor of the associated space $\mathbb A_N$:

  \begin{equation}
    R^i_{jmn}=L^i_{\underline{jm},n}-L^i_{\underline{jn},m}
    +
    L^\alpha_{\underline{jm}}L^i_{\underline{\alpha n}}-
    L^\alpha_{\underline{jn}}L^i_{\underline{\alpha m}}.
    \label{eq:R}
  \end{equation}

  Special symmetric affine connection spaces are Riemannian spaces
  $\mathbb R_N$ whose affine connection coefficients are the
  corresponding Christoffel symbols $\Gamma^i_{\underline{jk}}$ of
  the second kind.

  Many authors have developed the theory of symmetric affine
  connection spaces and mappings between them. Some of them
  are J. Mike\v s \cite{mik1, mik10, mik2, mik3, mik5, mik6, mik7,
  mik8}, N. S. Sinyukov \cite{sinjukov},
   V. E. Berezovski \cite{mik1, mik5, mik10, mik7, mik8}, L. P.
   Eisenhart \cite{eisRim} and many others. M. S. Stankovi\'c (see
   \cite{micasg}) obtained an invariant of an almost geodesic
   mapping of a non-symmetric affine connection space from the
   corresponding transformation of the curvature tensor
   (\ref{eq:R}).

    \subsection{Motivation}

    H. Weyl \cite{weyl} and T. Y. Thomas \cite{thomas1} obtained
    invariants of geodesic mappings of a symmetric affine connection
    space. J. Mike\v s with his
    research group \cite{mik1, mik10, mik2, mik3, mik5, mik6, mik7,
    mik8}, N. S. Sinyukov \cite{sinjukov}, and many other authors
    have continued the process of generalization of these
    invariants.

    The search for general formula for invariants of geometric
    mappings is started in \cite{vesicgeninv1}. The basic invariants
    of geometric mappings are obtained in this paper. It is also
    studied a special case in that paper. In this special case, it
    is obtained basic invariants of a mapping but the author founded
    some other invariants different of the basic ones. We are
    interested to obtain some of these invariants which are not
    basic in this paper.

    The invariants which we are interested to develop in this paper
    are

    \begin{itemize}
   \item The generalized Thomas projective parameter \cite{thomas1}:

  \begin{equation}
    T^i_{jk}=L^i_{\underline{jk}}-
    \frac1{N+1}\big(\delta^i_kL^\alpha_{\underline{j\alpha}}+
    \delta^i_jL^\alpha_{\underline{k\alpha}}\big).
    \label{eq:Thomasgeodesic}
  \end{equation}
  \item The Weyl projective tensor \cite{weyl}:

  \begin{equation}
    W^i_{jmn}=R^i_{jmn}+\frac1{N+1}\delta^i_jR_{[mn]}+
    \frac N{N^2-1}\delta^i_{[m}R_{jn]}+
    \frac1{N^2-1}\delta^i_{[m}R_{n]j}.
    \label{eq:Weylgeodesic}
  \end{equation}
  \end{itemize}

  In a Riemannian space $\mathbb{R}_N$, the Weyl
  projective tensor (\ref{eq:Weylgeodesic}) reduces to

  \begin{equation}
    W^i_{jmn}=R^i_{jmn}+
    \frac 1{N-1}\big(\delta^i_{m}R_{jn}-\delta^i_{n}R_{jm}\big).
    \label{eq:WeylgeodesicRN}
  \end{equation}

  \section{Reminder on basic invariants}

    Let $f:\mathbb{GA}_N\to\mathbb{G\overline A}_N$ be a
  mapping between non-symmetric affine connection spaces
  $\mathbb{GA}_N$ and $\mathbb{G\overline A}_N$. The deformation
  tensor $P^i_{jk}$ of
  this mapping satisfies the corresponding equation
  \cite{vesicgeninv1}

  \begin{equation}
    P^i_{jk}=\overline L^i_{jk}-L^i_{jk}=
    \overline\omega{}^i_{jk}-\omega^i_{jk}+\overline\tau{}^i_{jk}-
    \tau^i_{jk},
    \label{eq:P}
  \end{equation}

    \noindent for geometrical objects
  $\omega^i_{jk},\overline\omega{}^i_{jk},\tau^i_{jk},\overline\tau{}^i_{jk}$
  of the type $(1,2)$ such that $\omega^i_{jk}=\omega^i_{kj},$
  $\overline\omega{}^i_{jk}=\overline\omega{}^i_{kj},$
  $\tau^i_{jk}=-\tau^i_{kj},$
  $\overline\tau{}^i_{jk}=-\overline\tau{}^i_{kj}$.

    After symmetrizing
     the equation (\ref{eq:P}) by indices $j$ and
  $k$, one gets

  \begin{eqnarray}
    \overline L^i_{\underline{jk}}=L^i_{\underline{jk}}+
    \overline\omega^i_{jk}-\omega^i_{jk}.
    \label{eq:Psim}
  \end{eqnarray}

  With regard to the last result, it is obtained three kinds of
  invariants of the mapping $f$ in \cite{vesicgeninv1}.
  In this paper, we are interested to
  generalize just invariants of the second kind. The basic invariants
  which we will generalize are:

  \begin{align}
    &\aligned
    {\widetilde{\mathcal
    T}}{}^i_{jk}=L^i_{\underline{jk}}-\omega^i_{jk},
    \endaligned\label{eq:Tsimantisimgeneral}\\\displaybreak[0]
    &\aligned
    {\widetilde{\mathcal W}}{}^i_{jmn}=
    R^i_{jmn}-\omega{}^i_{jm|n}+
    \omega{}^i_{jn|m}+
    \omega{}^\alpha_{jm}\omega{}^i_{\alpha
    n}-\omega{}^\alpha_{jn}\omega{}^i_{\alpha
    m},
    \endaligned\label{eq:Wbasic}
  \end{align}

  \noindent for
  $p^1_1,\ldots,p^2_3=1,2,\omega_{(1).jk}^i=L^i_{\underline{jk}},
  \omega^i_{(2).jk}=\omega^i_{jk}$. The invariant
  (\ref{eq:Tsimantisimgeneral}) is the basic invariant of the
  mapping $f$ of the Thomas type but the invariant
  (\ref{eq:Wbasic}) is the basic invariant
  of the mapping $f$ of the Weyl type.




  \section{Derived invariants}

    In general, the geometrical object $\omega^i_{jk}$ from the
  equation (\ref{eq:P}) has the form

  \begin{equation}
    \omega^i_{jk}=s_{(1)}\big(\delta^i_j\rho_k+\delta^i_k\rho_j\big)+
    s_{(2)}\big(F^i_j\sigma_k+F^i_k\sigma_j\big)+s_{(3)}\sigma_{jk}\varphi^i,
    \label{eq:omegagen}
  \end{equation}

  \noindent for $s_1,s_2,s_3\in\mathbb R$, $1$-forms $\rho_j,\sigma_j$, an affinor $F^i_j$, a
  covariant tensor $\sigma_{jk}$ symmetric by indices $j$ and $k$
  and a contra-variant vector $\varphi^i$.

  \begin{rem}
    The geometrical objects $\delta^i_j\rho_k+\delta^i_k\rho_j$,
    $F^i_j\sigma_k+F^i_k\sigma_j$, $\sigma_{jk}\varphi^i$ are
    linearly independent. Otherwise, there would not be necessary
    three constants $s_{(1)},s_{(2)},s_{(3)}$.
  \end{rem}

  With regard to the equation (\ref{eq:omegagen}), we get

  \begin{equation}
    \aligned
    \omega^\alpha_{jm}\omega^i_{\alpha
    n}&=\big(s_{(1)}\big)^2\delta^i_j\rho_m\rho_n+\big(s_{(1)}\big)^2
    \delta^i_m\rho_j\rho_n\\&+
    \delta^i_n\Big(2\big(s_{(1)}\big)^2\rho_j\rho_m+
    s_{(1)}s_{(2)}\big(F^\alpha_m\rho_\alpha\sigma_j+F^\alpha_j\rho_\alpha\sigma_m\big)
    +s_{(1)}s_{(3)}\sigma_{jm}\rho_\alpha\varphi^\alpha\Big)\\&+
    \big(s_{(2)}\big)^2\Big(F^i_n\big(F^\alpha_m\sigma_j+F^\alpha_j\sigma_m\big)\sigma_\alpha+
    F^i_\alpha\big(F^\alpha_m\sigma_j+F^\alpha_j\sigma_m\big)\sigma_n\Big)+
    \big(s_{(3)}\big)^2\sigma_{jm}\sigma_{\alpha
    n}\varphi^\alpha\varphi^i\\&+
    s_{(1)}s_{(2)}\Big(F^i_n\big(\rho_j\sigma_m+\rho_m\sigma_j\big)+
    F^i_m\big(\rho_j\sigma_n+\rho_n\sigma_j\big)+
    F^i_j\big(\rho_m\sigma_n+\rho_n\sigma_m\big)\Big)\\&+
    s_{(1)}s_{(3)}\big(\sigma_{mn}\rho_j+\sigma_{jn}\rho_m+\sigma_{jm}\rho_n\big)\varphi^i\\&+
    s_{(2)}s_{(3)}\Big(\big(F^\alpha_m\sigma_j+F^\alpha_j\sigma_m\big)
    \sigma_{\alpha n}\varphi^i+\big(F^i_n\sigma_\alpha+
    F^i_\alpha\sigma_n\big)\sigma_{jm}\varphi^\alpha\Big)
    \endaligned\label{eq:omegaomegagen}
  \end{equation}

  It holds the following theorem:

  \begin{thm}
    Let $f:\mathbb{GA}_N\to\mathbb{G\overline A}_N$ be a geometric
    mapping. The geometrical objects

    \begin{align}
      &\aligned
      \widetilde{\mathcal T}{}^{(s).i}_{jk}=L^i_{\underline{jk}}-
      s_{(1)}\big(\delta^i_j\rho_k+\delta^i_k\rho_j\big)-
      s_{(2)}\big(F^i_j\sigma_k+F^i_k\sigma_j\big)-s_{(3)}\sigma_{jk}\varphi^i,
      \endaligned\label{eq:basicThomasGeneral}\\\displaybreak[0]
      &\aligned
      \widetilde{\mathcal
      W}{}^{(s).i}_{jmn}&=R^i_{jmn}-\delta^i_j\zeta^{(s)}_{[mn]}-
      \delta^i_m\zeta^{(s)}_{jn}+\delta^i_n\zeta^{(s)}_{jm}
      +\widetilde{\mathcal D}{}^{(s_2).(s_3).i}_{jmn}-\widetilde{\mathcal D}{}^{(s_2).(s_3).i}_{jnm},
      \endaligned\label{eq:basicWeylGeneral}
    \end{align}

    \noindent for

    \begin{align}
      &\aligned
      \zeta^{(s)}_{ij}&=s_{(1)}\rho_{i|j}+\big(s_{(1)}\big)^2\rho_i\rho_j+
      s_{(1)}s_{(2)}\big(F^\alpha_i\sigma_j+F^\alpha_j\sigma_i\big)\rho_\alpha+
      s_{(1)}s_{(3)}\sigma_{ij}\rho_\alpha\varphi^\alpha,
      \endaligned\label{eq:zeta(s)}\\
      &\aligned
      \widetilde{\mathcal
      D}{}^{(s_2).(s_3).i}_{jmn}&=\big(s_{(2)}\big)^2\Big(F^i_n\big(F^\alpha_m\sigma_j+
      F^\alpha_j\sigma_m\big)\sigma_\alpha+F^i_\alpha
      F^\alpha_m\sigma_j\sigma_n\Big)+\big(s_{(3)}\big)^2\sigma_{jm}\sigma_{\alpha
      n}\varphi^\alpha\varphi^i\\&+
      s_{(2)}s_{(3)}\Big(\big(F^\alpha_{m}\sigma_j+
      F^\alpha_j\sigma_{m}\big)\sigma_{\alpha n}\varphi^i-
      \big(F^i_{m}\sigma_\alpha+
      F^i_\alpha\sigma_{m}\big)\sigma_{jn}\varphi^\alpha\Big)\\&-
      s_{(2)}\big(F^i_j\sigma_m+F^i_m\sigma_j\big)_{|n}-
      s_{(3)}\big(\sigma_{jm}\varphi^i\big)_{|n},
      \endaligned\label{eq:D(s2)(s3)}
    \end{align}

    \noindent $s=(s_{1},s_{2},s_{3})$, are the basic invariants of the mapping $f$.\hfill\qed
  \end{thm}

  \subsection{Invariants in associated space}

  We will analyze the invariants (\ref{eq:basicThomasGeneral},
  \ref{eq:basicWeylGeneral}) of a mapping $f:\mathbb{GA}_N\to\mathbb{G\overline A}_N$
  bellow. From this analyzing, we will
  obtain some other invariants of the mapping $f$.

  After contracting the equality $\widetilde{\overline{\mathcal T}}{}^i_{jk}-
  \widetilde{\mathcal T}{}^i_{jk}=0$ by indices $i$ and $k$, we get

  \begin{equation}
    \aligned
    (N+1)s_{(1)}\big(\overline\rho_j-\rho_j\big)&=\overline
    L^\alpha_{\underline{j\alpha}}-s_{(2)}\big(\overline
    F^\alpha_j\overline\sigma_k+\overline F\overline\sigma_j\big)-
    s_{(3)}\overline\sigma_{j\alpha}\overline\varphi^\alpha\\&-
    L^\alpha_{\underline{j\alpha}}+s_{(2)}\big(F^\alpha_j\sigma_\alpha+F\sigma_j\big)+
    s_{(3)}\sigma_{j\alpha}\varphi^\alpha,
    \endaligned\label{eq:rho-rho}
  \end{equation}

  \noindent for $F=F^\alpha_\alpha,\overline F=\overline
  F^\alpha_\alpha$.

  After substituting this equation into the equality $\widetilde{\overline{\mathcal T}}{}^{(s).i}_{jk}-
  \widetilde{\mathcal T}{}^{(s).i}_{jk}=0$, one obtains

  \begin{equation*}
    \widetilde{{\overline T}}{}^{(s).i}_{jk}=
    \widetilde{T}{}^{(s).i}_{jk},
  \end{equation*}

  \noindent for

  \begin{equation}
    \aligned
    \widetilde T{}^{(s).i}_{jk}&=L^i_{\underline{jk}}-
    \frac{s_{(1)}}{N+1}\delta^i_j\Big(L^\alpha_{\underline{k\alpha}}-
    s_{(2)}(F^\alpha_k\sigma_\alpha+F\sigma_k\big)-s_{(3)}\sigma_{k\alpha}\varphi^\alpha\Big)\\&
    -
    \frac{s_{(1)}}{N+1}\delta^i_k\Big(L^\alpha_{\underline{j\alpha}}-
    s_{(2)}(F^\alpha_j\sigma_\alpha+F\sigma_j\big)-s_{(3)}\sigma_{j\alpha}\varphi^\alpha\Big)\\&-
    s_{(2)}\big(F^i_j\sigma_k+F^i_k\sigma_j\big)-s_{(3)}\sigma_{jk}\varphi^i,
    \endaligned\label{eq:generalThomassymmetric}
  \end{equation}

  \noindent and the corresponding $\widetilde{\overline
  T}{}^{(s).i}_{jk}$.

  \begin{lem}
    Let $f:\mathbb{GA}_N\to\mathbb{G\overline A}_N$ be a geometric
    mapping of a non-symmetric affine connection space
    $\mathbb{GA}_N$. The geometrical object
    \emph{(\ref{eq:generalThomassymmetric})} is an invariant of the
    mapping $f$.\hfill\qed
  \end{lem}

  \begin{cor}
    The invariant \emph{(\ref{eq:generalThomassymmetric})} and the
    Thomas projective parameter $T^i_{jk}$ given in the equation
    \emph{(\ref{eq:Thomasgeodesic})} satisfy the equation

    \begin{equation}
      \aligned
          \widetilde T{}^{(s).i}_{jk}&=s_{(1)}
          T^{i}_{jk}+(1-s_{(1)})L^i_{\underline{jk}}
          -
    s_{(2)}\big(F^i_j\sigma_k+F^i_k\sigma_j\big)-s_{(3)}\sigma_{jk}\varphi^i
          \\&+
    \frac{s_{(1)}}{N+1}\Big(
    s_{(2)}(F^\alpha_k\sigma_\alpha+F\sigma_k\big)+s_{(3)}\sigma_{k\alpha}\varphi^\alpha\Big)\delta^i_j\\&
    +
    \frac{s_{(1)}}{N+1}\Big(
    s_{(2)}(F^\alpha_j\sigma_\alpha+F\sigma_j\big)+s_{(3)}\sigma_{j\alpha}\varphi^\alpha\Big)\delta^i_k,
      \endaligned\label{eq:dThomasThomascorrelation}
    \end{equation}

    \noindent for the corresponding $s=(s_1,s_2,s_3)$.\hfill\qed
  \end{cor}

  The geometrical object (\ref{eq:generalThomassymmetric}) is
  \emph{the derived associated invariant} of the Thomas
  type.

  Furthermore, from the invariance of the geometrical object (\ref{eq:basicWeylGeneral}) we
  obtain

  \begin{equation}
    \aligned
    \overline
    R^i_{jmn}&=R^i_{jmn}+\delta^i_j\big(\overline\zeta^{(s)}_{[mn]}-
    \zeta^{(s)}_{[mn]}\big)+\delta^i_m\big(\overline\zeta^{(s)}_{jn}-
    \zeta^{(s)}_{jn}\big)-\delta^i_n\big(\overline\zeta^{(s)}_{jm}-
    \zeta^{(s)}_{jm}\big)\\&-
    \widetilde{\overline{\mathcal D}}{}^{(s_2).(s_3).i}_{jmn}+
    \widetilde{\overline{\mathcal D}}{}^{(s_2).(s_3).i}_{jnm}+
    \widetilde{{\mathcal D}}{}^{(s_2).(s_3).i}_{jmn}-
    \widetilde{{\mathcal D}}{}^{(s_2).(s_3).i}_{jnm}.
    \endaligned\label{eq:R-Rgeneral}
  \end{equation}

  After contracting this equation by indices $i$ and $j$, we get

  \begin{equation}
    \aligned
    (N+1)\big(\overline\zeta^{(s)}_{[mn]}-\zeta^{(s)}_{[mn]}\big)&=
    -\overline R_{[mn]}+\widetilde{\overline{\mathcal
    D}}{}^{(s_2).(s_3).\alpha}_{\alpha mn}-
    \widetilde{\overline{\mathcal
    D}}{}^{(s_2).(s_3).\alpha}_{\alpha nm}\\&+R_{[mn]}-
    \widetilde{{\mathcal
    D}}{}^{(s_2).(s_3).\alpha}_{\alpha mn}+
    \widetilde{{\mathcal
    D}}{}^{(s_2).(s_3).\alpha}_{\alpha nm}.
    \endaligned\label{eq:R-Rgenerali=j}
  \end{equation}

  From the equations (\ref{eq:R-Rgeneral}, \ref{eq:R-Rgenerali=j}),
  one obtains

  \begin{equation}
    \aligned
    \overline
    R^i_{jmn}&=R^i_{jmn}+\frac1{N+1}\delta^i_j\big(R_{[mn]}-\overline
    R_{[mn]}\big)+\delta^i_m\big(\overline\zeta^{(s)}_{jn}-
    \zeta^{(s)}_{jn}\big)-\delta^i_n\big(\overline\zeta^{(s)}_{jm}-
    \zeta^{(s)}_{jm}\big)\\&
    +\frac1{N+1}\delta^i_j\big(\widetilde{\overline{\mathcal
    D}}{}^{(s_2).(s_3).\alpha}_{\alpha mn}-
    \widetilde{\overline{\mathcal
    D}}{}^{(s_2).(s_3).\alpha}_{\alpha nm}-
    \widetilde{{\mathcal
    D}}{}^{(s_2).(s_3).\alpha}_{\alpha mn}+
    \widetilde{{\mathcal
    D}}{}^{(s_2).(s_3).\alpha}_{\alpha nm}\big)\\&-
    \widetilde{\overline{\mathcal D}}{}^{(s_2).(s_3).i}_{jmn}+
    \widetilde{\overline{\mathcal D}}{}^{(s_2).(s_3).i}_{jnm}+
    \widetilde{{\mathcal D}}{}^{(s_2).(s_3).i}_{jmn}-
    \widetilde{{\mathcal D}}{}^{(s_2).(s_3).i}_{jnm}.
    \endaligned\label{eq:R-Rgenerali=j2}
  \end{equation}

  From the contraction of this result by indices $i$ and $n$, we get

  \begin{equation}
    \aligned
    (N-1)\big(\overline\zeta^{(s)}_{jm}-\zeta^{(s)}_{jm}\big)&=
    \frac N{N+1}\big(R_{jm}-\overline R_{jm}\big)+
    \frac1{N+1}\big(R_{mj}-\overline R_{mj}\big)\\&+
    \frac1{N+1}\big(
    \widetilde{\overline{\mathcal
    D}}{}^{(s_2).(s_3).\alpha}_{\alpha[jm]}-
    \widetilde{{\mathcal
    D}}{}^{(s_2).(s_3).\alpha}_{\alpha[jm]}\big)\\&-
    \widetilde{\overline{\mathcal D}}{}^{(s_2).(s_3).\alpha}_{j[m\alpha]}+
    \widetilde{{\mathcal D}}{}^{(s_2).(s_3).\alpha}_{j[m\alpha]}.
    \endaligned\label{eq:R-Rgenerali=n}
  \end{equation}

  With regard to the equations (\ref{eq:R-Rgenerali=j2},
  \ref{eq:R-Rgenerali=n}), we get

  \begin{equation*}
    {{\overline W}}{}^{(s).[1].i}_{jmn}=
    {{W}}{}^{(s).[1].i}_{jmn},
  \end{equation*}

  \noindent for

  \begin{equation}
    \aligned
    {{W}}{}^{(s).[1].i}_{jmn}&=R^i_{jmn}+
    \frac1{N+1}\delta^i_jR_{[mn]}+
    \frac N{N^2-1}\delta^i_{[m}R_{jn]}+
    \frac1{N^2-1}\delta^i_{[m}R_{n]j}\\&
    +\widetilde{\mathcal D}{}^{(s_2).(s_3).i}_{j[mn]}-
    \frac1{N+1}\delta^i_j\widetilde{\mathcal D}{}^{(s_2).(s_3).\alpha}_{\alpha[mn]}\\&+
    \frac1{N^2-1}\delta^i_m\big((N+1)\widetilde{{\mathcal
    D}}{}^{(s_2).(s_3).\alpha}_{j[n\alpha]}-\widetilde{{\mathcal
    D}}{}^{(s_2).(s_3).\alpha}_{\alpha[jn]}\big)\\&-
    \frac1{N^2-1}\delta^i_n\big((N+1)\widetilde{{\mathcal
    D}}{}^{(s_2).(s_3).\alpha}_{j[m\alpha]}-\widetilde{{\mathcal
    D}}{}^{(s_2).(s_3).\alpha}_{\alpha[jm]}\big),
    \endaligned\label{eq:generalWeylderived1}
    \tag{$i$}
  \end{equation}

  \noindent and the corresponding ${{\overline
  W}}{}^{(s).[1].i}_{jmn}$.

  Let us test are some summands in the invariant
  (\ref{eq:generalWeylderived1}) invariants of the mapping $f$.
  After contract the equality ${{\overline
  W}}{}^{(s).[1].i}_{jmn}-{{W}}{}^{(s).[1].i}_{jmn}=0$
  by indices $i$ and $n$, we get

  \begin{equation*}
    \widetilde{\overline{\mathcal
    D}}{}^{(s_2).(s_3).\alpha}_{\alpha[jm]}-
    \widetilde{{\mathcal
    D}}{}^{(s_2).(s_3).\alpha}_{\alpha[jm]}=0.
  \end{equation*}

  Hence, the invariant (\ref{eq:generalWeylderived1}) reduces to

    \begin{equation}
    \aligned
    {{W}}{}^{(s).[2].i}_{jmn}&=R^i_{jmn}+
    \frac1{N+1}\delta^i_jR_{[mn]}+
    \frac N{N^2-1}\delta^i_{[m}R_{jn]}+
    \frac1{N^2-1}\delta^i_{[m}R_{n]j}\\&
    +\widetilde{\mathcal D}{}^{(s_2).(s_3).i}_{j[mn]}+
    \frac1{N-1}\big(\delta^i_m\widetilde{{\mathcal
    D}}{}^{(s_2).(s_3).\alpha}_{j[n\alpha]}-
    \delta^i_n\widetilde{{\mathcal
    D}}{}^{(s_2).(s_3).\alpha}_{j[m\alpha]}\big),
    \endaligned\label{eq:generalWeylderived1final}
    \tag{$ii$}
  \end{equation}

\noindent  and the corresponding
${\overline{W}}{}^{(s).[2].i}_{jmn}$.

  Let us check-out are there invariants of the mapping $f$ into the invariant
  (\ref{eq:generalWeylderived1final}). After contracting the
  equality ${\overline{W}}{}^{(s).[2].i}_{jmn}-
  {{W}}{}^{(s).[2].i}_{jmn}=0$ by the indices $i$ and $j$,
  one obtains

  \begin{equation*}
    \widetilde{\overline{\mathcal
    D}}{}^{(s_2).(s_3).\alpha}_{j[n\alpha]}-
    \widetilde{{\mathcal
    D}}{}^{(s_2).(s_3).\alpha}_{j[n\alpha]}=0.
  \end{equation*}

  Hereof, the invariant (\ref{eq:generalWeylderived1final}) reduces
  to

  \begin{equation}
  \aligned
    W{}^{(s).i}_{jmn}&=R^i_{jmn}+
    \frac1{N+1}\delta^i_jR_{[mn]}+
    \frac N{N^2-1}\delta^i_{[m}R_{jn]}+
    \frac1{N^2-1}\delta^i_{[m}R_{n]j}\\&
    +\widetilde{\mathcal D}{}^{(s_2).(s_3).i}_{jmn}-
    \widetilde{\mathcal D}{}^{(s_2).(s_3).i}_{jnm}.
  \endaligned\label{eq:Weylgensimf}
  \end{equation}

  It holds the following theorem:

  \begin{thm}
    Let $f:\mathbb{GA}_N\to\mathbb{G\overline A}_N$ be a geometric
    mapping. The geometrical object\linebreak
    \emph{(\ref{eq:Weylgensimf})} is an invariant of this mapping.
    \hfill\qed
  \end{thm}

  \begin{cor}
    The invariant \emph{(\ref{eq:Weylgensimf})} and the Weyl
    projective tensor \emph{(\ref{eq:Weylgeodesic})} satisfy the
    equation

    \begin{equation}
      W{}^{(s).i}_{jmn}=W^i_{jmn}
    +\widetilde{\mathcal D}{}^{(s_2).(s_3).i}_{jmn}-
    \widetilde{\mathcal D}{}^{(s_2).(s_3).i}_{jnm},
    \label{eq:dWeylWeylcorrelation}
    \end{equation}

    \noindent for the corresponding $s=(s_1,s_2,s_3)$.\hfill\qed
  \end{cor}

  The geometrical object (\ref{eq:Weylgensimf}) is \emph{the derived
  associated invariant} of the Weyl type.

  \subsection{$F$-planar mappings}

  In this part of paper, we will apply the above obtained results.
  The theoretical part of this application will be search for invariants
  of conformal mappings. The practical example will be about a
  transformation of a special Riemannian space.\\[3pt]

  \noindent\textbf{$F$-planar mappings.} A mapping $f:\mathbb{A}_N\to\mathbb{\overline A}_N$ is called
  \emph{the $F$-planar mapping} if it is determined with the
  equation

  \begin{equation}
    \overline L^i_{\underline{jk}}=
    L^i_{\underline{jk}}+\delta^i_k\psi_j+\delta^i_j\psi_k+
    F^i_k\sigma_j+F^i_j\sigma_k,
    \label{eq:F-planardfn}
  \end{equation}

  \noindent for $1$-forms $\psi_j,\sigma_j$ and an affinor $F^i_j$.

  The basic equation of the inverse mapping $f^{-1}:\mathbb
  {A}_N\to\mathbb{\overline{A}}_N$ is

  \begin{equation}
    L^i_{\underline{jk}}=\overline L^i_{\underline{jk}}-
    \delta^i_k\psi_j-\delta^i_j\psi_k-F^i_k\sigma_j-
    F^i_j\sigma_k.
    \label{eq:F-planar-1dfn}
  \end{equation}

  Hence, we get the mapping $f^{-1}$ is an $F$-planar mapping for

  \begin{eqnarray}
    \overline F^i_j=F^i_j,&
    \overline\sigma_j=-\sigma_j,&
    \overline\psi_j=-\psi_j.
    \label{eq:F-planarpsisigmaF-1}
  \end{eqnarray}

  Moreover, it holds

  \begin{equation}
    \overline L^i_{\underline{jk}}=
    L^i_{\underline{jk}}+\delta^i_k\psi_j+\delta^i_j\psi_k-
    \frac12\overline F^i_k\overline\sigma_j-
    \frac12\overline F^i_j\overline\sigma_k+
    \frac12F^i_k\sigma_j+\frac12F^i_j\sigma_k.
    \label{eq:F-planaromega-omega1}
  \end{equation}

  Therefore, the $F$-planar mapping $f$ is the case of
  $s_1=1,s_2=1/2,s_3=0$,
  in the equation (\ref{eq:omegagen}).

  After contracting the equation (\ref{eq:F-planaromega-omega1}) by
  indices $i$ and $k$, one gets

  \begin{equation}
  \aligned
    \psi_j&=\frac1{N+1}\Big(\overline L^\alpha_{\underline{j\alpha}}
    +\frac12\overline F^i_k\overline\sigma_j+
    \frac12\overline F^i_j\overline\sigma_k\Big)-
    \frac1{N+1}\Big(L^\alpha_{\underline{j\alpha}}
    +\frac12F^i_k\sigma_j+
    \frac12F^i_j\sigma_k\Big).
  \endaligned\label{eq:F-planarpsij}
  \end{equation}

  With regard to the equations (\ref{eq:omegagen},
  \ref{eq:F-planaromega-omega1},
  \ref{eq:F-planarpsij}), we obtain

  \begin{equation}
    \rho_j=\frac1{N+1}\Big(L^\alpha_{\underline{j\alpha}}
    +\frac12F\sigma_j+
    \frac12F^\alpha_j\sigma_\alpha\Big),
    \label{eq:F-planarrhoj}
  \end{equation}

  \noindent for $F=F^\alpha_\alpha$.

  From the equation (\ref{eq:basicThomasGeneral}), it holds that the
  associated invariant of the Thomas type of the mapping $f$ is

  \begin{align}
  &\aligned
    \widetilde{T}{}^i_{jk}&=L^i_{\underline{jk}}-\frac12\big(F^i_j\sigma_k+F^i_k\sigma_j\big)-
    \frac1{N+1}\delta^i_j\Big(L^\alpha_{\underline{k\alpha}}-
    \frac12\big(F^\alpha_k\sigma_\alpha+F\sigma_k\big)\Big)
    \\&-
    \frac1{N+1}\delta^i_k\Big(L^\alpha_{\underline{j\alpha}}-
    \frac12\big(F^\alpha_j\sigma_\alpha+F\sigma_j\big)\Big),
  \endaligned\label{eq:F-planarThomas}\\&
  \aligned
    \widetilde{T}{}^i_{jk}&=T^i_{{jk}}-\frac12\big(F^i_j\sigma_k+F^i_k\sigma_j\big)+
    \frac1{2(N+1)}\delta^i_j
    \big(F^\alpha_k\sigma_\alpha+F\sigma_k\big)
    \\&+
    \frac1{2(N+1)}\big(F^\alpha_j\sigma_\alpha+F\sigma_j\big),
  \endaligned\tag{\ref{eq:F-planarThomas}'}\label{eq:F-planarThomas'}
  \end{align}

  \noindent for the generalized Thomas projective parameter
  $T^i_{jk}$ from the equation (\ref{eq:Thomasgeodesic}).

  Based on the equation (\ref{eq:F-planarpsisigmaF-1}), we obtain

  \begin{equation*}
    \overline F^i_j\overline
    F^m_n\overline\sigma_p\overline\sigma_q=
    F^i_jF^m_n\sigma_p\sigma_q,
  \end{equation*}

  \noindent i.e. the geometrical objects (\ref{eq:zeta(s)}, \ref{eq:D(s2)(s3)}) reduce
  to

  \begin{align}
    &\aligned\widetilde{\mathcal
    D}{}^i_{jmn}=-\frac12\big(F^i_j\sigma_m+F^i_m\sigma_j\big)_{|n},
    \endaligned\label{eq:F-planarD(s2)(s3)}\\
    &\aligned
    \zeta_{ij}&=\frac1{N+1}L^\alpha_{\underline{i\alpha}|j}+\frac1{(N+1)^2}
    L^\alpha_{\underline{i\alpha}}L^\beta_{\underline{j\beta}}+
    \frac1{2(N+1)}L^\beta_{\underline{\alpha\beta}}
    \big(F^\alpha_i\sigma_j+F^\alpha_j\sigma_i\big)\\&+
    \frac1{2(N+1)}\Big(\big(F\sigma_i+F^\alpha_i\sigma_\alpha\big)_{|j}+
    L^\alpha_{\underline{i\alpha}}\big(F\sigma_j+F^\beta_j\sigma_\beta\big)+
    L^\alpha_{\underline{j\alpha}}\big(F\sigma_i+F^\beta_i\sigma_\beta\big)\Big).
    \endaligned\label{eq:F-planarsigma}
  \end{align}

  Therefore, the invariant (\ref{eq:basicWeylGeneral}) of the
  mapping $f$ is

  \begin{equation}
    \aligned
    \widetilde{\mathcal
    W}{}^i_{jmn}&=R^i_{jmn}+\frac1{N+1}\delta^i_jR_{[mn]}-\frac12\big(F^i_j\sigma_m+F^i_m\sigma_j\big)_{|n}+
    \frac12\big(F^i_j\sigma_n+F^i_n\sigma_j\big)_{|m}-\delta^i_{[m}\zeta_{jn]},
    \endaligned\label{eq:F-planarWbasic}
  \end{equation}

  \noindent for $\zeta_{ij}$ given in the equation
  (\ref{eq:F-planarsigma}).

  The invariant (\ref{eq:Weylgensimf}) of the $F$-planar mapping $f$
  is

  \begin{equation}
    \widetilde W^i_{jmn}=W^i_{jmn}-\frac12\big(F^i_j\sigma_m+F^i_m\sigma_j\big)_{|n}+
    \frac12\big(F^i_j\sigma_n+F^i_n\sigma_j\big)_{|m},
    \label{eq:F-planarWderived}
  \end{equation}

  \noindent for the Weyl projective tensor $W^i_{jmn}$.\\[3pt]

  \begin{exm} In this example, we are aimed to obtain the invariants
  \emph{(\ref{eq:F-planarThomas}, \ref{eq:F-planarWbasic},
  \ref{eq:F-planarWderived})}
  of an $F$-planar mapping $f:\mathbb R_3\to\mathbb{\overline R}_3$ for
  the  Riemannian space $\mathbb R_3$ determined with the metric
  tensor

  \begin{equation}
    g_{\underline{ij}}=\left[\begin{array}{ccc}
      u^2&0&0\\
      0&v^2&0\\
      0&0&w^2
    \end{array}\right]
  \end{equation}

  The corresponding affinor $F^i_j$ and covariant vector
  $\sigma_j$ are

  \begin{eqnarray}
    F^i_j=\left[\begin{array}{ccc}
      \sin u&0&0\\
      0&\cos v&0\\
      0&0&w
    \end{array}\right]&\mbox{and}&
    \sigma_j=\left[\begin{array}{c}
      0\\
      0\\
      \ln(1+u^2+v^2+w^2)
    \end{array}\right]
  \end{eqnarray}

  Let be $\mathcal F^i_{jk}=F^i_k\sigma_j+F^i_j\sigma_k$.
  It is satisfied

  \begin{equation}
    \mathcal F^i_{jk}=\left\{\begin{array}{ll}
      0,&j,k\in\{1,2\}\\
      \sin u\ln(1+u^2+v^2+w^2),&i=j=1,k=3\mbox{ or }i=k=1,j=3,\\
      \cos v\ln(1+u^2+v^2+w^2),&i=j=2,k=3\mbox{ or }i=k=2,j=3,\\
      2w\ln(1+u^2+v^2+w^2),&i=j=k=3.
    \end{array}\right.
    \label{eq:F-planarFijk}
  \end{equation}

  It also holds

  \begin{equation}
    \mathcal F^\alpha_{j\alpha}=\big(\sin u+\cos
    v+w)\sigma_j+F^{(j)}_j\sigma_{(j)}
    \label{eq:F-planarFj}
  \end{equation}

    The Christoffel symbols of the second kind of this space
  are

  \begin{equation}
    \begin{array}{ccc}
      \Gamma^1_{\underline{11}}=\dfrac1u,&\Gamma^1_{\underline{22}}=\dfrac
      v{u^2},&\Gamma^1_{\underline{33}}=\dfrac w{u^2},\\\\
      \Gamma^2_{\underline{11}}=\dfrac u{v^2},&
      \Gamma^2_{\underline{22}}=\dfrac1v,&
      \Gamma^2_{\underline{33}}=\dfrac w{v^2},\\\\
      \Gamma^3_{\underline{11}}=\dfrac u{w^2},&
      \Gamma^3_{\underline{22}}=\dfrac v{w^2},&
      \Gamma^3_{\underline{33}}=\dfrac 1w,
    \end{array}
  \end{equation}

  \noindent and $\Gamma^i_{\underline{jk}}=0$ in all other cases.

  The generalized Thomas projective parameter of the space $\mathbb R_3$ is

  \begin{equation}
    T^i_{jk}=\left\{\begin{array}{ll}
      \Gamma^i_{\underline{jk}},&i\not\in\{j,k\},\\
      -\frac1{N+1}\Gamma^{(k)}_{\underline{k(k)}},&i=j\neq k,\\
      \frac{N-1}{N+1}\big(\frac1u\delta^i_1+\frac1v\delta^i_2+\frac1w\delta^i_3\big),&i=j=k.
    \end{array}\right.
  \end{equation}

   Hence, the derived invariant of Thomas type of the mapping $f$ is

  \begin{equation}
    \widetilde{\mathcal T}{}^i_{jk}=T^i_{jk}-\frac12\mathcal
    F^i_{jk}+\frac1{8}\delta^i_j\mathcal F^\alpha_{k\alpha}+
    \frac18\delta^i_k\mathcal F^\alpha_{j\alpha},
  \end{equation}

  \noindent for $\mathcal F^i_{jk},\mathcal F^\alpha_{j\alpha}$ given by the
  equations
  \emph{(\ref{eq:F-planarFijk}, \ref{eq:F-planarFj})}.

  We have the following cases for the curvature
  tensor $R^i_{jmn}$:

  \begin{equation}
    \aligned
    1.&\quad m=n\Rightarrow R^i_{jmn}=0\Rightarrow R_{jm}=0,\\
    2.&\quad j=m\neq n\Rightarrow R^i_{jmn}=
    \Gamma^i_{\underline{jm},n}+\Gamma^{(n)}_{\underline{jm}}\Gamma^i_{\underline{(n)n}}
    \Rightarrow R_{jm}={\Gamma^{\alpha}_{\underline{jm},\alpha}}+
    \Gamma^{\alpha}_{\underline{jm}}L^{(\alpha)}_{\underline{\alpha(\alpha)}},\\
    3.&\quad j=n\neq m\Rightarrow
    R^i_{jmn}=-\Gamma^i_{\underline{jn},m}-\Gamma^{\alpha}_{\underline{jn}}
    \Gamma^i_{\underline{\alpha m}}\Rightarrow R_{jm}=-
    \Gamma^{(m)}_{\underline{j(j)}}\Gamma^{(j)}_{\underline{(m)m}},\\
    4.&\quad j\neq m,j\neq n\Rightarrow R^i_{jmn}=0\Rightarrow
    R_{jm}=0.
    \endaligned
  \end{equation}

  It is also satisfied

  \begin{equation}
    \aligned
    \zeta_{ij}&=\frac1{N+1}\Gamma^{(i)}_{\underline{i(i)}|j}+
    \frac1{(N+1)^2}\Gamma^{(i)}_{\underline{i(i)}}
    \Gamma^{(j)}_{\underline{j(j)}}\\&+\frac1{2(N+1)}
    \Big(\Gamma^{(\alpha)}_{\underline{\alpha(\alpha)}}
    \mathcal F^\alpha_{ij}+
    \mathcal F^\alpha_{i\alpha|j}+
    \Gamma^{(i)}_{\underline{i(i)}}\mathcal F^\alpha_{j\alpha}+
    \Gamma^{(j)}_{\underline{j(j)}}
    \mathcal F^\alpha_{i\alpha}\Big),
    \endaligned\label{eq:F-planarexamplezeta}
  \end{equation}

  \noindent for $\mathcal F^i_{jk},\mathcal F^\alpha_{j\alpha}$
  given in the equations \emph{(\ref{eq:F-planarFijk},
  \ref{eq:F-planarFj})}.

  Hence, the corresponding invariants of Weyl type of the mapping $f$ are

  \begin{align}
  &\aligned
  \widetilde{\mathcal W}{}^i_{jmn}=R^i_{jmn}-\delta^i_{[m}\zeta_{jn]}-\frac12
  \mathcal F^i_{jm|n}+\frac12\mathcal F^i_{jn|m},
  \endaligned\\
  &\aligned
    \widetilde{
    W}{}^i_{jmn}&=R^i_{jmn}+\frac1{N-1}\big(\delta^i_mR_{jn}-\delta^i_nR_{jm}\big)
    -\frac12\mathcal F^i_{jm|n}+
    \frac12\mathcal F^i_{jn|m}.
  \endaligned
  \end{align}
\end{exm}


\begin{thebibliography}{33}

      \bibitem{mik7} \textbf{V. Berezovski, S. B\'asc\'o, J. Mike\v
    s}, \emph{Almost geodesic mappings of affinely connected spaces
    that preserve the Riemannian curvature}, Annales Mathematicae et
    Informaticae, 45 (2015) pp. 3--10.

    \bibitem{mik10} \textbf{V. Berezovski, S. B\'asc\'o, J. Mike\v
    s}, \emph{Diffeomorphism of Affine Connected Spaces Which
    Preserved Riemannian and Ricci Curvature Tensors}, Miskolc
    Mathematical Notes, Vol. 18 (2017), No. 1, pp. 117--124.

    \bibitem{mik8} \textbf{V. Berezovskij, J. Mike\v s}, \emph{On
    special almost geodesic mappings of type $\pi_1$ of spaces with
    affine connection}, Acta Universitatis Palackianae Olomucensis.
    Facultas Rerum Naturalium. Mathematica, Vol. 43 (2004), No.1,
    21--26.

    \bibitem{eisRim} \textbf{L. P. Eisenhart}, \emph{Riemannian
    Geometry}, Princeton University Press, 1926.

         \bibitem{e1}
  {\bf A. Einstein}, \emph{A generalization of the
  relativistic theory of gravitation}, Ann. of. Math., 45 (1945),
  No. 2, 576--584.

  \bibitem{e2}
  {\bf A. Einstein}, \emph{Bianchi identities in the
  generalized theory of gravitation}, Can. J. Math., (1950), No. 2,
  120--128.

  \bibitem{e3}
  {\bf A. Einstein}, \emph{Relativistic Theory of
  the Non-symmetric Field}, Princeton University Press, New Jersey,
  1954, 5th edition.

            \bibitem{mik3} \textbf{J. Mike\v s}, \emph{Holomorphically
    Projective Mappings and Their Generalizations}, International
    Journal of Mathematical Sciences, Vol. 89, No. 3, 1998,
    1334--1353.

    \bibitem{mik1} \textbf{J. Mike\v s, V. E. Berezovski, E.
    Stepanova, H. Chud\'a}, \emph{Geodesic Mappings and Their
    Generalizations}, Journal of Mathematical Sciences, Vol. 217.,
    No. 5, 2016, 607--623.

    \bibitem{mik6} \textbf{J. Mike\v s, V. Kiosak, A. Van\v
    zurov\'a}, \emph{Geodesic mappings of manifolds with affine
    connection}, Olomouc: Palacky University, 2008.

        \bibitem{mik5} \textbf{J. Mike\v s, E. Stepanova,
    A. Van\v zurova, et al.}, \emph{Differential geometry of special mappings},
    Olomouc: Palacky University, 2015.

      \bibitem{mik2} \textbf{J. Mike\v s, A. Van\v zurov\'a, I.
    Hinterleitner}, \emph{Geodesic mappings and some
    generalizations}, Olomouc: Palacky University, 2009.



       \bibitem{sinjukov} \textbf{N. S. Sinyukov}, \emph{Geodesic mappings
    of Riemannian spaces},  (in Russian), "Nauka", Moscow, 1979.

    \bibitem{micasg} {\textbf{M. S. Stankovi\'c}, \emph{Special equitorsion
almost geodesic mappings of the third type of non-symmetric affine
connection spaces}, Appl. Math. and Computation, 244 (2014),
595--701.}

          \bibitem{thomas1} \textbf{T. Y. Thomas}, \emph{On the
    projective and equi-projective geometries of paths},
    Proc. Nat. Acad. Sc. 11 (1925) 199--203.

    \bibitem{vesicgeninv1} \textbf{N. O. Vesi\'c}, \emph{Invariants of Geometric Mappings},
    ArXiv: 2422024.

            \bibitem{weyl} \textbf{H. Weyl}, \emph{Zur infinitesimal
    geometrie: Einordnung der projectiven und der konformen
    auffssung}, Gottingen Nachrichten (1921), 99--112.
\end{thebibliography}
\end{document}